\documentclass{amsart}

\usepackage{pb-diagram, amssymb}
\usepackage[english]{babel}
\newtheorem{theorem}{Theorem}[section]
\newtheorem*{theorem*}{Theorem}
\newtheorem{lemma}[theorem]{Lemma}
\newtheorem{proposition}[theorem]{Proposition}

\theoremstyle{definition}
\newtheorem{definition}[theorem]{Definition}
\newtheorem{example}[theorem]{Example}
\newtheorem{examples}[theorem]{Examples}
\newtheorem{question}[theorem]{Question}

\newcommand{\proven}{\qedhere}

\newcommand{\id}{\ensuremath{\mathrm{id}}}
\newcommand{\supp}{\ensuremath{\mathop{\mathrm{supp}}}}

\newcommand{\abs}[1]{\left| #1 \right|}
\newcommand{\norm}[1]{\left\| #1 \right\|}
\newcommand{\isom}{\cong}

\newcommand{\holalg}[1]{\mathcal{O}({#1})}
\newcommand{\holalgb}[1]{\mathcal{B}({#1})}
\newcommand{\enmo}[1]{\ensuremath{{\mathrm{End}\!\left( {#1} \right)}}}
\newcommand{\aut}[1]{\ensuremath{{\mathrm{Aut}\left( {#1} \right)}}}

\newcommand{\CC}{{\mathbb{C}}}
\newcommand{\RR}{{\mathbb{R}}}

\newcommand{\ZZ}{{\mathbb{Z}}}
\newcommand{\NN}{{\mathbb{N}}}
\newcommand{\PP}{{\mathbb{P}}}
\newcommand{\EE}{{\mathbb{E}}}

\begin{document}
\title{Stein Spaces Characterized by their Endomorphisms}
\author{Rafael B. Andrist}
\date{2nd November 2008}
\address{Mathematisches Institut, Universit\"at Bern, Sidlerstr. 5, 3012 Bern, Switzerland}
\email{rafael.andrist@math.unibe.ch}
\thanks{This work was partially supported by Schweizerischer Nationalfonds grant no. 200021-116165}
\begin{abstract}
Finite dimensional Stein spaces admitting a proper holomorphic embedding of the complex line are characterized, among all complex spaces, by their holomorphic endomorphism semigroup in the sense that any semigroup isomorphism induces either a biholomorphic or an antibiholomorphic map between them.
\end{abstract}
\maketitle

\section{Introduction}
The question if an object is determined by the algebraic structures ``naturally'' arising on it, such as the rings/algebras of functions, is quite old. A theorem of Bers \cite{Bers} (for planar domains), Nakai \cite{Nakai} (for Riemann surfaces with non-constant bounded functions) and Iss'sa \cite{Isssa} (for normal reduced Stein spaces) states that if the algebras of holomorphic functions $\holalg{X}, \holalg{Y}$ are isomorphic as abstract rings, then there exists a unique biholomorphic or antibiholomorphic map $\varphi : X \to Y$ inducing this isomorphism by $f \mapsto f \circ \varphi^{-1}$; if in these cases they are isomorphic as algebras, the map $\varphi$ is biholomorphic. 
\par
Other algebraic structures which belong to an object, are its automorphism group and its semigroup of endomorphisms. The question of characterizing a complex manifold by its automorphism group was studied for some special cases, with the additional assumption however, that the automorphism groups are not only isomorphic as (abstract) groups, but also isomorphic as topological groups. This has been done by Isaev for the unit ball \cite{Isaev1}, \cite{Isaev2}, for the polydisc \cite{Isaev3} and by Isaev and Kruzhilin for $\CC^n$ \cite{Isaev4}. Without topology, the situation is much more complicated, and it may happen that some domains cannot be distinguished anymore. This leads to the study of the semigroup of endomorphisms of a complex space.
\par
It is easy to see that a biholomorphic or antibiholomorphic $\varphi : X \to Y$ induces an isomorphism of the holomorphic endomorphism semigroups, by conjugating $f : X \to X$ to $\varphi \circ f \circ \varphi^{-1} : Y \to Y$. On the other hand, each semigroup isomorphism determines a unique map $\varphi : X \to Y$ which induces the isomorphism by conjugating the endomorphisms (see proposition \ref{prop_conjug}). It is therefore natural to ask if all such isomorphisms are necessarly induced by a biholomorphic or antibiholomorphic $\varphi : X \to Y$.
\par
For topological spaces and continous endomorphisms, the analogon to this question has been studied in detail, see e.g. the survey by Magill \cite{Magill}. But for complex manifolds only a few results are known.
Hinkkanen \cite{Hinkkanen} showed this in 1992 for $\varphi : \CC \to \CC$ and gave at the same time counterexamples for certain unbounded domains in $\CC$, which may have a semigroup consisting only of the identity and the constant maps. In 1993, Eremenko \cite{Eremenko} proved this for Riemann surfaces admitting non-constant bounded functions, generating ``enough'' endomorphisms with certain properties using the one-dimensional Schr\"oder equation. His result was generalized in 2002 for bounded domains in $\Omega_1, \Omega_2 \subset \CC^n$ by Merenkov \cite{Merenkov}. He first noted that $\varphi : \Omega_1 \to \Omega_2$ is a homeomorphism, then he locally linearized the map $\varphi$ by finally reducing this problem to a one-dimensional Schr\"oder equation again, and cleverly excluded the possibilities that $\varphi$ could be holomorphic in one direction and anti-holomorphic in another. The same proof was applied to $\varphi : \CC^n \to Y, \; n \geq 2,$ by Buzzard and Merenkov \cite{Buzzard2}, the additional difficulty was to show that $\varphi$ is a homeomorphism. This was possible thanks to the generation of a sub-basis of the topology through certain Fatou-Bieberbach domains constructed before by Buzzard and Hubbard \cite{Buzzard1}. 
\par
We observe that all complex manifolds for which the result has been established, are Stein spaces or at least their function algebra is Stein (in case of domains in $\CC^n$). In a more general setting, generalizations of the proof of Eremenko are limited by the problem to show that the map $\varphi$ is a homeomorphism. Therefore we want to make use of the theory of Stein spaces, and focus more on analytic subsets instead of showing first that $\varphi$ is a homeomorphism. The already mentioned theorem of Bers, Nakai and Iss'sa, which states that if the algebras of holomorphic functions $\holalg{X}, \holalg{Y}$ are isomorphic as abstract rings, then there exists a unique biholomorphic or antibiholomorphic map $\varphi : X \to Y$ inducing this isomorphism by $f \mapsto f \circ \varphi^{-1}$, motivates the ``emulation'' of holomorphic functions by holomorphic endomorphisms. However, our method of proof does not need the theorem of Iss'sa, as our situation is simplified by the fact that the map $\varphi : X \to Y$ is already given by the semigroup isomorphism. Our main result is
\begin{theorem*}[\ref{thm_conjug1}]
Let $X$ and $Y$ be complex spaces and $\varphi: X \to Y$ such that $f \mapsto \varphi \circ f \circ \varphi^{-1}$ is an isomorphism between the endomorphism semigroups of $X$ and $Y$. Then $\varphi$ is either biholomorphic or antibiholomorphic if the following criteria are fulfilled:
\begin{enumerate}
\item $X$ is a finite dimensional Stein space with countable topology
\item $X$ admits a proper holomorphic embedding $i : \CC \hookrightarrow X$.
\end{enumerate}
\end{theorem*}
In particular, this theorem covers the case of $\CC^n$, but also a lot more classes of examples, including the known examples of manifolds with density property, such as certain homogenous spaces with actions of a semi-simple Lie group, linear algebraic spaces and hypersurfaces of the form $\{(x, u, v) \in \CC^n \times \CC \times \CC \,:\, f(x) = u \cdot v\} $. The result can be slightly generalized, as done in theorem \ref{thm_conjug2}, and works with some modifications also for embeddings of the unit disk $\EE$, which is dealt with in theorem \ref{thm_conjug3}.
\par
Throughout the the paper, we will from now on assume that all complex spaces have a countable topology without mentioning it explicitly.

\newpage
\section{Endomorphisms and Conjugation}

\begin{definition}
Let $X$ be a complex space. The set of holomorphic endomorphisms $f : X \rightarrow X$, equipped with the composition of functions as multiplication, is a semigroup with the identity as neutral element and it is called $\enmo{X}$. Similarly, for base-point preserving endomorphisms $f : (X, x_0) \to (X, x_0)$ it is called $\enmo{X, x_0}$. We also introduce the notation for the automorphism group $\aut{X}$, and $\aut{X, x_0}$ for the base-point preserving automorphisms.
\end{definition}

\par The following proposition is due to Buzzard and Merenkov \cite{Buzzard2}, but except for the part about weak double-transitivity, it was already known in 1937, see Schreier \cite{Schreier}.

\begin{proposition}
\label{prop_conjug}
Let $X$ and $Y$ be complex spaces. Let $\Phi : \enmo{X} \rightarrow \enmo{Y}$ be a homomorphism of semigroups.
Then there exists a bijective map $\varphi : X \rightarrow Y$ such that $$\Phi(f) = \varphi \circ f \circ \varphi^{-1} \quad \forall \, f \in \enmo{X}$$
if one of the following conditions is fulfilled:
\begin{itemize}
\item $\Phi$ is an isomorphism of semigroups.
\item $\Phi$ is an epimorphism of semigroups and $X$ is weakly double-transitive, i.e. for any two different $x_1, x_2 \in X$ and $x_1^\prime \in X$ there exists an open set $U \subset X$ with $x_1^\prime \in U$ such that for all $x_2^\prime \in U \setminus x_1^\prime$ there exists $f \in \enmo{X}$ with $f(x_1) = x_1^\prime, \; f(x_2) = x_2^\prime$.
\end{itemize}
\end{proposition}
\begin{proof}
First note that a map $f:X \rightarrow X$ is constant if and only if $f \circ g = f \; \forall\, g \in \enmo{X}$.
Therefore, an epimorphism $\Phi : \enmo{X} \rightarrow \enmo{Y}$ maps constant maps to constant maps. 
This induces in a natural way a map $\varphi : X \rightarrow Y$ by $\varphi(x) = y$ where $y$ is such that $c_y = \Phi(c_x)$
and $c_x$ resp. $c_y$ denote the constant maps with image point $x$ resp. $y$.

\par Because $\Phi$ is an epimorphism of semigroups and because of the one-to-one relation between points and constant maps, $\varphi$ is surjective. If $\Phi$ were an isomorphism, $\varphi$ would automatically become a bijective map. However, using the weak double-transitivity one can show that $\varphi$ is injective without this additional condition on $\Phi$: Choose $x_1 \in X$ such that in any neighbourhood $U$ of $x_1$ there is always a point $x_U^\prime$ with $\varphi(x_1) \neq \varphi(x^\prime)$. Assume by contradiction that there is more than one point in the pre-image of $y = \varphi(x_1) \in Y$, i.e. $x_1, x_2 \in \varphi^{-1}(y), x_1 \neq x_2$. There must be an endomorphism $f : X \to X$ with $f(x_1) = x_1$ and $f(x_2) = x_U^\prime \not \in \varphi^{-1}(y)$. Using constant maps, this results in
$$
\begin{array}{l}
\Phi(f \circ c_{x_1}) = \Phi(f) \circ \Phi(c_{x_1}) = \Phi(f) \circ c_{y} = \Phi(c_{x_1}) \\
\Phi(f \circ c_{x_2}) = \Phi(f) \circ \Phi(c_{x_2}) = \Phi(f) \circ c_{y} = \Phi(c_{x^\prime}) \neq \Phi(c_{x_1})
\end{array}
$$

\par Finally, look at the maps $\varphi^{-1} \circ \Phi(f) \circ \varphi : X \rightarrow X$:
$$
\begin{array}{clclcl}
             & \varphi^{-1} \circ \Phi(f) \circ \varphi \circ c_x
         & = & \varphi^{-1} \circ \Phi(f) \circ \Phi(c_x) 
         & = & \varphi^{-1} \circ \Phi(f \circ c_x) \\
           = & \varphi^{-1} \circ \Phi(c_{f(x)})
         & = & \varphi^{-1} \circ c_{\varphi(f(c_x))}
         & = & c_{f(x)} \quad \forall \, x \in X \\
\end{array}
$$
This is equivalent to $\Phi(f) = \varphi \circ f \circ \varphi^{-1}$. \proven
\end{proof}

\par This motivates the following

\begin{definition}
A (set-theoretic) map $\varphi : X \to Y$ between two complex spaces is called a \textit{conjugating map} if it is bijective and if it induces a homomorphism $\Phi : \enmo{X} \to \enmo{Y}$ of the endomorphism semigroups by conjugating $\enmo{X} \ni f \mapsto \varphi \circ f \circ \varphi^{-1} \in \enmo{Y}$. If, in addition, $\Phi$ is an isomorphism, then $\varphi$ is called an \textit{iso-conjugating map}.
\end{definition}

\par From now on, we will only talk about (iso-)conjugating maps and not consider the semigroup isomorphisms in an abstract way anymore.

\section{Conjugation for Stein Spaces with a complex line}

\begin{proposition}
\label{prop_conjonedim}
Let $\varphi : \CC \to \CC$ be a conjugating map. Then $\varphi \in \langle \aut{\CC}, z \mapsto \overline{z} \rangle$.
\end{proposition}
\begin{proof}
By composing $\varphi$ with an automorphism (i.e. a non-degenerate affine $\CC$-linear map) one may assume that $\varphi(0) = 0$ and $\varphi(1) = 1$.
\begin{enumerate}
\item $\varphi$ is a field automorphism: Since $\varphi : \CC \to \CC$ is a conjugating map, it maps automorphisms to automorphisms, i.e.
$$\varphi(a \cdot w + b) = A(a, b) \cdot \varphi(w) + B(a, b), \; w = \varphi^{-1}(z)$$
Setting $w = 0$ we get $B(a, b) = \varphi(b)$, independent of $a$.
Setting $w = 1$ and $b = 0$ we get $\varphi(a) = A(a, 0)$, thus $\varphi(a \cdot w) = \varphi(a) \cdot \varphi(w)$, meaning that $\varphi$ is multiplicative and $\varphi(-1) = -1$.
Finally we choose $a \cdot w + b = 0$, i.e. $w = -b/a$. Then $A(a,b) \cdot \varphi(-b/a) + \varphi(b) = 0$ which resolves to $A(a, b) = \varphi(a)$ and implies additivity of $\varphi$.
\item The field isomorphism $\varphi$ is continous: It is sufficient to prove continuity in $0$. Assume by contradiction that there is a sequence $(c_j)_{j \in \NN}$ with $\lim_{j \to \infty} c_j = 0$, but $\lim_{j \to \infty} \varphi(c_j) \neq 0$. Set $f(z) := \sum_{j \in \NN} (c_j)^j z^j$. This is an entire function and must be mapped to an entire function $\varphi \circ f \circ \varphi^{-1}(z) = \sum_{j \in \NN} (\varphi(c_j))^j z^j$. But this requires $\lim_{j \to \infty} \varphi(c_j) = 0$. Contradiction.
\end{enumerate}
A continous field automorphism of $\CC$ is either the identity or the complex conjugation. \proven
\end{proof}

\par This has already been proven by Hinkkanen \cite{Hinkkanen}, but with a more complicated proof, because he did not made complete use of the explicit form of the automorphism group. Note that it is not required that $\varphi$ is an iso-conjugating map which will be important in the proof of theorem \ref{thm_conjug2}. Also note that if $\varphi$ would be assumed to be continous, it would be sufficient to look at the automorphism group alone.

\begin{lemma}
\label{lem_analyticsub}
Let $X$ be a finite dimensional Stein space with a holomorphic embedding $i : \CC \hookrightarrow X$. Then for a subset $A \subseteq X$ the following is equivalent:
\begin{enumerate}
\item $A$ is an analytic set in $X$.
\item $\exists \,x_0 \in X, \; \exists\, F_1, \dots F_q \in \enmo{X}$ such that
\[A = \left\{x \in X \,:\, F_1(x) = \dots = F_q(x) = x_0\right\}\]
\end{enumerate}
The implication $2 \Longrightarrow 1$ holds for all complex spaces.
\end{lemma}
\begin{proof} \hfill
\begin{itemize}
\item[] $1 \Longrightarrow 2$: Because $X$ is a Stein space, there are finitely many holomorphic functions $f_k : X \to \CC, \; k = 1, \dots, q$ such that
\[
A = \left\{x \in X \,:\, f_1(x) = \dots = f_q(x) = 0\right\}
\]
As shown by Forster and Ramspott \cite{ForsterRamspott}, if $X$ has complex dimension $n \in \NN$, then $n$ holomorphic functions are sufficient. Set $x_0 := i(0)$. Define $F_k : X \to X$ by $F_k = i \circ f_k, \; k = 1, \dots, q$. \\
\item[] $2 \Longrightarrow 1$: $F_k^{-1}(x_0) \subseteq X, \; k = 1, \dots, q$ is obviously closed and is locally the finite intersection of zeros of holomorphic functions $f_{1,k}, \dots f_{q_k,k} : U \subseteq X \to \CC$. Therefore, $F_k^{-1}(x_0) \subset X$ is an analytic subset, and
\[
A = \left\{x \in X \;:\; f_{\ell_k,k}(x) = 0, \; \ell_k = 1, \dots, q_k, \; k = 1, \dots q \right\}
\]
is an analytic subset too. \proven
\end{itemize}
\end{proof}

\begin{theorem}
\label{thm_conjug1}
Let $X$ and $Y$ be complex spaces and $\varphi: X \to Y$ an iso-conjugating map. Then $\varphi$ is either biholomorphic or antibiholomorphic if the following criteria are fulfilled:
\begin{enumerate}
\item $X$ is a finite dimensional Stein space
\item $X$ admits a proper holomorphic embedding $i : \CC \hookrightarrow X$.
\end{enumerate}
\end{theorem}

\par The result is a consequence of the following, slightly more general, theorem:

\begin{theorem}
\label{thm_conjug2}
Let $X$ and $Y$ be complex spaces and $\varphi: X \to Y$ an iso-conjugating map. Then $\varphi$ is either biholomorphic or antibiholomorphic if the following criteria are fulfilled:
\begin{enumerate}
\item There is a holomorphic (not necessarly proper) embedding $X \hookrightarrow \CC^n$.
\item $X$ admits a holomorphic embedding $i : \CC \hookrightarrow X$ such that all holomorphic endomorphisms of $i(\CC) \isom \CC$ are restrictions of holomorphic endomorphisms of $X$ and such that $i(\CC) = \bigcap_{k = 1, \dots, q} G_k^{-1}(x_0)$ for a point $x_0 \in X$ and $G_1, \dots, G_q \in \enmo{X}$.
\item There is a holomorphic function $\omega : X \to \CC$ which is non-constant on $i(\CC)$.
\end{enumerate}
\end{theorem}
\begin{proof} \hfill
\begin{enumerate}
\item The image $B := \varphi(i(\CC))$ can be written as 
\[
\varphi(i(\CC)) = \bigcap_{k= 1, \dots, q} \left(\varphi \circ G_k \circ \varphi^{-1} \right)^{-1}(\varphi(x_0)),
\]
therefore it is an analytic subset of $Y$ (by lemma \ref{lem_analyticsub}), and carries the structure of a complex space.
\item We want to show that $\varphi|i(\CC) \to B$ is a conjugating map: By assumption the endomorphisms of the analytic subset $i(\CC) \isom \CC$ can be extended to holomorphic endomorphisms $F : X \to X$ and are conjugated to endomorphisms of $Y$ which restrict to endomorphisms of $B$. The goal is to apply proposition \ref{prop_conjonedim}, but for this we need to show that $B \isom \CC$ first.
\item The group $\aut{\CC}$ acts double-transitively on $\CC$ and so does $\aut{B}$ on $B$. Therefore the structure around each point in $B$ is locally the same and $B$ is in fact a complex manifold. Next we show that $B$ is connected:
\[
\begin{diagram}
\node[2]{X} \arrow{e,t}{\varphi} \arrow{sw,t}{\omega} \arrow{s,r}{} \node{Y} \arrow{s,r}{\varphi \circ i \circ \omega \circ \varphi^{-1}} \\
\node{\CC} \arrow{e,t}{i} \arrow{s,r}{} \node{i(\CC)} \arrow{e,t}{\varphi} \arrow{s,l}{} \node{B} \arrow{s,r}{G} \\
\node{\CC} \arrow{e,t}{i} \node{i(\CC)} \arrow{e,t}{\varphi} \node{B}
\end{diagram}
\]
Assume by contradiction that there are several connected components. Then there is a holomorphic map $G : B \to B$ which is the identity on one component, but constant on all other components. The map $i \circ \omega : X \to i(\CC) \subseteq X$ is non-constant when restricted to $i(\CC)$, and for its conjugate we have $\varphi \circ i \circ \omega \circ \varphi^{-1} : Y \to B \subseteq Y$. The map $G \circ \varphi \circ i \circ \omega \circ \varphi^{-1} : Y \to B$ is holomorphic on $Y$ and still non-constant when restricted to $B$. Conjugating it back and restricting to $i(\CC)$, we end up with a non-constant map $\CC \to \CC$ which omits an infinite number of points. This is a contradiction to Picard's theorem and we conclude that $B$ is connected.
\item $(\CC, +)$ is an abelian complex Lie group. We show that $B$ also carries the structure of an abelian complex Lie group by defining an addition $\oplus: B \times B \to B$ 
\[
a \oplus b = \varphi\left( \varphi^{-1}(a) + \varphi^{-1}(b) \right)
\]
This operation obviously gives $(B, \oplus)$ an abelian group structure. For fixed $a \in B$, the map $b \mapsto a \oplus b$ is conjugated to a translation in $\CC$ and therefore holomorphic, similarly $a \mapsto a \oplus b$ is holomorphic too. This is sufficient to make $(B, \oplus)$ a complex Lie group.\footnote{Separate continuity is already sufficient for a topological group (Montgomery \cite{Montgomery}), provided the underlying metric space is countable and locally complete. Then Hartogs theorem for separate analyticity can be used. Note that without continuity, separate analyticity is not enough for complex manifolds. See Palais \cite{Palais} for a brief overview.}
A connected abelian complex Lie group is isomorphic to a product $\CC^m \times {\CC^*}^n \times T$ with $\holalg{T} = \CC$, and $T$ is trivial if and only if $X$ is Stein (see e.g. Grauert and Remmert \cite{GrauertRemmert}, and Matsushima and Marimoto \cite{Matsushima2}).
Both ${\CC^*}^n$ and $T$ contain, if non-trivial, a lot of elements of finite order, but $\CC^m$ does not. By construction, $\varphi|i(\CC)$ is a group isomorphism (not a priori a Lie group isomorphism) and therefore $B \isom \CC^m$. The case $B \isom \CC^m, m \geq 2,$ can be ruled out in the same way as the non-connectedness: Choose $G : B \to B$ as a projection to a coordinate axis in $\CC^m \isom B$. Thus, $B \isom \CC$, the proposition \ref{prop_conjonedim} can be applied, and as a result, $\varphi|i(\CC)$ is either biholomorphic or anti-biholomorphic to its image.
\item The algebras $\holalg{X}$ and $\holalg{Y}$ of holomorphic functions are isomorphic as $\CC$-algebras by sending $\holalg{X} \ni f \mapsto (\varphi \circ i)^{-1} \circ (\varphi \circ i \circ f \circ \varphi^{-1}) = f \circ \varphi^{-1} =: \hat{\Phi}(f) \in \holalg{Y}$, if necessary with complex conjugation. The properties of a $\CC$-algebra homomorphism are fulfilled, since addition and multiplication of functions is point-wise. $\hat{\Phi}$ is an isomorphism because the situation in $X$ and $Y$ is symmetric and because $\varphi$ is an iso-conjugating map. In case that an additional complex conjugation is necessary for the holomorphic functions, we replace $X$ by its image under complex conjugation in $\CC^n$.
\[
\begin{diagram}
\node[2]{X} \arrow{e,t}{\varphi} \arrow{sw,t}{f} \arrow{s,r}{} \node{Y} \arrow{s,l}{} \arrow{se,t}{f \circ \varphi^{-1}} \\
\node{\CC} \arrow{e,t}{i} \node{X} \arrow{e,t}{\varphi} \node{Y} \node{\CC} \arrow{w,t}{\varphi \circ i}
\end{diagram}
\]
\item By assumption, there is a holomorphic embedding $(f_1, \dots, f_n) : X \hookrightarrow \CC^n$ with $f_1, \dots, f_n \in \holalg{X}$. It follows that $f_k \circ \varphi^{-1} \in \holalg{Y}, \; k=1,\dots,n$ and that $(f_1 \circ \varphi^{-1}, \dots, f_n \circ \varphi^{-1}) : Y \hookrightarrow  \CC^n$ is a holomorphic embedding too (but not neccesarly proper). As the map $(f_1, \dots, f_n)$ is biholomorphic to its image, this implies that $\varphi^{-1}$ is holomorphic. Now we apply the same argument to $\hat{\Phi}^{-1} : \holalg{Y} \to \holalg{X}$ and obtain that $\varphi$ is holomorphic too. \proven
\end{enumerate}
\end{proof}
\par The second condition can be reformulated in the language of sheaf theory: There is a coherent ideal $\mathfrak{I} \subseteq \mathcal{O}_X$ which is generated by its global sections and such that $\supp(\mathcal{O}_X / \mathfrak{I}) \isom \CC$ and $H^1(X, \mathfrak{I}) = 0$.
\begin{proof}[Proof of theorem \ref{thm_conjug1}] \hfill \\
The complex space $X$ is Stein and therefore there is a proper holomorphic embedding $X \hookrightarrow \CC^n$ (see e.g. Grauert and Remmert \cite{GrauertRemmert}). By assumption, $i : \CC \hookrightarrow X$ is a proper holomorphic embedding. Then $i(\CC)$ is a closed subset of $X$ and locally the zero set of finitely many holomorphic functions and therefore an analytic subset of $X$ -- and since $X$ is Stein, we can apply lemma \ref{lem_analyticsub} to find the desired endomorphisms $G_1, \dots, G_q \in \enmo{X}$.  Finally, an endomorphism of $f : i(\CC) \to i(\CC)$ can be viewed as a holomorphic function $i^{-1} \circ f : i(\CC) \to \CC$. This is the restriction of some function $F : X \to \CC$, which is a direct consequence of Cartan's Theorem B, because $i(\CC)$ is an analytic subset of a Stein space (see e.g. Grauert and Remmert \cite{GrauertRemmert} again). Then $i \circ F \in \enmo{X}$ has the property that $i \circ F|i(\CC) = f$. The existence of a holomorphic map $\omega : X \to \CC$, non-constant on $i(\CC)$, follows directly from the definition of a Stein space, as there are functions separating points.\proven
\end{proof}

\begin{example}[for the neccessity of the conditions] \hfill \\
Here is an example (mainly due to Hinkkanen \cite{Hinkkanen}) which shows that $\varphi : X \to Y$ being an iso-conjugating map and a homeomorphism between Stein manifolds is not sufficient to guarantee (anti-)holomorphicity: Set $X := \PP^1 \setminus \{ 2^n, n \in \NN \}, Y := \PP^1 \setminus \{ 3^n, n \in \NN \}$.
Endomorphisms of $X$ and $Y$ can be continued to endomorphisms of $\PP^1$: No essential singularities are possible, otherwise the image could omit at most two points in $\PP^1$ by the theorem of Picard. The sequences of ``gaps'' in $\PP^1$ are exponential, but possible endomorphisms are only polynomial, therefore $\enmo{X} = \{ \id_X \} \cup X$ and $\enmo{Y} = \{ \id_Y \} \cup Y$, and in addition there is no biholomorphic or antibiholomorphic map $X \to Y$, which would need to be the restriction of a (possibly complex conjugated) Moebius transformation $\PP^1 \to \PP^1$.
But then any bijective map $\varphi : X \to Y$ is an iso-conjugating map and it is easy to choose it as a homeomorphism, e.g. piece-wise linear for strips $[2^n, 2^{n+1}] \times i \RR \mapsto [3^n, 3^{n+1}] \times i \RR$.
\end{example}

\newpage

\begin{examples}[for theorem \ref{thm_conjug1}]
\label{ex_conjug1} \hfill
\begin{enumerate}
\item \label{ex_cn} $\CC^n, \; n \in \NN$
\item $\CC \times X$, where $X$ is any Stein space: The embedding $i : \CC \to \CC \times X, \; z \mapsto (z, x_0), \; x_0 \in X,$ is obviously proper.
\item \label{ex_cstar} $\CC^* \times \CC^*$: The embedding $i : \CC \to \CC^* \times \CC^*, \; z \mapsto (e^{z}, e^{- \sqrt{-1} z})$ is proper: Note that the logarithms of the moduli of the component functions correspond to the real and imaginary part of $z$.
\item $X_1 \times X_2$, where $X_1$ and $X_2$ are Stein spaces and both have a proper embedding of $\CC^*$: $i_k : \CC^* \hookrightarrow X_k, \; k = 1, 2$. Then $i : \CC \hookrightarrow X_1 \times X_2, \; z \mapsto (i_1(e^{z}),i_2(e^{-\sqrt{-1}z}))$.
\item \label{ex_slinear} $\mathrm{SL}_n(\CC) \subset \CC^{n^2}$: The space is Stein as an analytic subset and e.g. $i : \CC \to \mathrm{SL}_n(\CC), \; z \mapsto E + z \cdot E_{1,n}$ is a proper holomorphic embedding.
\item \label{ex_semisimple} Homogenous Stein spaces $X = G/K$, where $G$ is a semi-simple Lie group. It is known that for a holomorphic action of a semi-simple (even for a reductive) Lie group $G$ with $X = G/K$ Stein, the subgroup $K$ is reductive and $X$ is affine algebraic; further in this algebraic realization the action of $G$ is algebraic (see Matsushima \cite{Matsushima1}, \cite{Matsushima2}; Borel and Harish-Chandra \cite{Borel}). There is an algebraic subgroup $H < G$ isomorphic to $\mathrm{SL}_2(\CC)$ and $H$ contains a properly embedded $\CC$ as in example \ref{ex_slinear}.
\item \label{ex_alggroup} Linear algebraic groups except $\CC^*$: There is a decomposition as a semi-direct product due to Mostow \cite{Mostow} such that $X = R \rtimes U$, where $R$ is reductive and $U$ is unipotent. $U$ is biholomorphic to $\CC^m$ which already gives the desired embedding if $m \neq 0$. $R$ contains a maximal semi-simple subgroup (previous example) or is biholomorphic to ${\CC^*}^m$. In the latter case, the $\CC$ is obtained for $m \geq 2$ as in example \ref{ex_cstar}.
\item \label{ex_hypersurface} A hypersurface $H \subset \CC^{n+2}$ of the form $H = \{ (x, u, v) \in \CC^n \times \CC \times \CC \,:\, f(x) = u \cdot v^m \}, \; n, m \in \NN,$ where $f : \CC^n \to \CC$ is a holomorphic function. Two cases need to be considered: If $f$ has a zero in $x_0 \in \CC$, then $i : \CC \hookrightarrow H, \; z \mapsto (x_0, z, 0)$ is a proper holomorphic embedding. If $f$ has no zeros at all, then there is even a $\CC^n$ properly embedded: $\CC^n \ni z \mapsto (z, f(z), 1) \in H$.
The complex space $H$ is a manifold if and only if $f'(x) \neq 0$ whenever $f(x) = 0$.
\end{enumerate}
\end{examples}
Note that the examples \ref{ex_semisimple} ($G$ with trivial center), \ref{ex_alggroup} (except ${\CC^*}^n$) and \ref{ex_hypersurface} (with $m=1$) are all known examples of Stein manifolds with density property (see section \ref{sec_questions_density} for more details).

\begin{examples}[for theorem \ref{thm_conjug2}] \hfill
\begin{enumerate}
\item $\CC \times \Omega \subseteq \CC^{m+1}$, where $\Omega \subset \CC^m$ is any open subset. The embedding $i : \CC \hookrightarrow \CC \times \Omega, \; z \mapsto (z, x_0), \; x_0 \in \Omega,$ is proper holomorphic, $i(\CC) = f^{-1}(0, x_0)$ for an endomorphism $f(z, w) = (0, w)$ and any $\alpha : i(\CC) \to i(\CC)$ is the restriction of an endomorphism $A : \CC \times \Omega \to \CC \times \Omega, \; A = \alpha \times \id_{\Omega}$.
\item $\CC^n \! \setminus \! \{ 0 \} \subset \CC^n, \; n \geq 2$: It is for sure not a domain of holomorphy and therefore not Stein. But the embedding $i : \CC \hookrightarrow \CC^n \! \setminus \! \{ 0 \}, \; z \mapsto (z, 1, \dots, 1)$ is proper, and for the endomorphism $f(z_1, \dots, z_n) = (1, z_2,\dots,z_n)$ we have $i(\CC) = f^{-1}(1,\dots,1)$, and any $\alpha : i(\CC) \to i(\CC)$ is the restriction of an endomorphism $A : \CC^n \! \setminus \! \{ 0 \} \to \CC^n \! \setminus \! \{ 0 \}, \; A(z, w) = (\alpha(z),1\dots,1)$.
\end{enumerate}
\end{examples}

\section{Conjugation for bounded domains}

For bounded domains in $\CC^n$, Merenkov \cite{Merenkov} proved the following result:
\begin{theorem}
\label{thm_merenkov}
Let $\Omega_1 \subseteq \CC^{n_1}, \Omega_2 \subseteq \CC^{n_2}$ be bounded domains. Then an iso-conjugating map $\varphi : \Omega_1 \to \Omega_2$ is either biholomorphic or antibiholomorphic, and $n_1 = n_2$.
\end{theorem}
\par
We first note that that this result can be reformulated in a stronger way: We do not need to require $\Omega_2$ to be a bounded domain:
\begin{theorem}
Let $\Omega \subseteq \CC^n$ be a bounded domain, and $Y$ a complex manifold. Then an iso-conjugating map $\varphi : \Omega \to Y$ is either biholomorphic or antibiholomorphic.
\end{theorem}
\begin{proof} It is sufficient to show that there is a biholomorphic map $Y \hookrightarrow \CC^n$ which realizes $Y$ as an open bounded subset.
\begin{enumerate}
\item By shifting $\Omega$ biholomorphically, we may assume $0 \in \Omega$.
\item $\varphi$ is a homeomorphism: For each point $x_0 \in \Omega$ and each $\delta > 0$ there is an $\varepsilon > 0$ such that $f : \Omega \to \Omega, \; x \mapsto \varepsilon \cdot x + x_0,$ is an injective holomorphic map whith an image contained in a ball of radius $\delta$ around $x_0$. Thus the Fatou-Bieberbach domains of $\Omega$ form a basis of the topology, and by lemma \ref{lem_homeo} (below), the map $\varphi : \Omega \to Y$ is a homeomorphism.
\item $A := \CC \times \{0\}^{n-1} \cap \Omega$ is an analytic subset of $\Omega$. There is an $a > 0$ such that $a \cdot \EE^n \subset \Omega$.
Define $F := a \cdot \left( 0, \frac{\pi_2}{\norm{\pi_2}_{\Omega}}, \dots, \frac{\pi_n}{\norm{\pi_n}_{\Omega}} \right) : \Omega \to \Omega$. Then $A = F^{-1}(0)$, and 
$B := \varphi(A) = \left(\varphi \circ F \circ \varphi^{-1}\right)^{-1}(\varphi(0))$. Therefore, $B \subset Y$ is an analytic subset (by lemma \ref{lem_analyticsub}) and in fact a complex space of dimension $1$, since $\varphi$ is a homeomorphism.
\item Possible singularities of $B$ are isolated points, and so their pre-images under $\varphi$ are isolated as well. Therefore we find a $b > 0$ such that $D := b \cdot \EE \times \{ 0 \}^{n-1}$ does not hit any of them. $\varphi(D)$ is then an open subset of a complex space, simply connected, without singularities and relatively compact. There must be a biholomorphic map $g : \varphi(D) \to \EE$.
\item Define $f_k : \Omega \to D \subset \Omega, \; k = 1, \dots, n,$ by $f_k := b \cdot \left(\frac{\pi_k}{\norm{\pi_k}_\Omega}, 0 \dots, 0\right)$. The maps $f_k$ ar conjugated to holomorphic maps $\varphi \circ f_k \circ \varphi^{-1} : Y \to \varphi(D)$, and $(g \circ \varphi \circ f_1 \circ \varphi^{-1}, \dots, g \circ \varphi \circ f_n \circ \varphi^{-1}) : Y \to \EE^n$ is the desired embedding. \proven
\end{enumerate}
\end{proof}
\par
In this section we will use methods similar to the ones from the previous section in order to show another way of proofing Merenkov's results for various bounded domains.

\begin{definition}
Let $X$ be a complex manifold. An open subset $\Omega \subset X$ is called a Fatou-Bieberbach domain (of $X$), if there exists a biholomorphic map $f : X \to \Omega$ and $\Omega \neq X$.
\end{definition}

\begin{lemma}
\label{lem_homeo}
Let $X, Y$ be complex manifolds (possibly infinite-dimensional, but with countable topology) and $\varphi : X \to Y$ a conjugating map. Assume the Fatou-Bieberbach domains of $X$ form a sub-basis of the topology. Then $\varphi$ is a homeomorphism.
\end{lemma}
\begin{proof} \hfill
\begin{enumerate}
\item Each Fatou-Bieberbach domain $\Omega \subset X$ is the image of an injective holomorphic map $f : X \to X$. By assumption any open set $U \subseteq X$ contains a finite intersection of such domains, i.e. $U \supseteq \bigcap_{j=1, \dots, n} f_j(X)$. These maps $f_j$ are conjugated to $g_j := \varphi \circ f_j \circ \varphi^{-1} : Y \to Y$ which are again injective holomorphic maps and therefore have an open image $g_j(Y) \subseteq Y$ and thus $\varphi(U) \supseteq \bigcap_{j=1, \dots, n} g_j(Y)$ where $\bigcap_{j=1, \dots, n} g_j(Y)$ is open.
The map $\varphi^{-1} : Y \to X$ is continous; for simpler notation, we switch the roles of $X$ and $Y$ and can assume that $\varphi : X \to Y$ is bijective and continous in the following purely topological arguments. It is sufficient to show $\varphi^{-1}$ is proper.
\item Let $\left(K_j\right)_{j \in \NN}$ be an exhaustion of $X$ by compacts $K_j \subset X$,
i.e. $K_j \subset {K_{j+1}}^\circ \; \forall \, j \in \NN$ and $\bigcup_{j \in \NN} K_j = X$.
Such an exhaustion exists because of local compactness and countable topology.
By $\varphi$, the $K_j$ are mapped to compacts $L_j := \varphi(K_j) \subset Y$ with $\bigcup_{j \in \NN} L_j = Y$ and $L_j \subset L_{j+1}$, but a priori not necessarly $L_j \subset L_{j+1}^\circ$.
\item The restricted map $\varphi|K_j$ is a homeomorphism onto its image. For all $x \in K_j^\circ$ there is an open neighborhood $U \subseteq K_j^\circ$ which can be considered as an open set in $\CC^n$, and similarly for $f(x) \in L_j$, there is an open neighborhood $V \subset Y$ which can be considered as an open set in $\CC^m$. Openness in the ambient euclidean space as well as the dimension is preserved by a homeomorpism, therefore $\varphi( U \cap \varphi^{-1}(V) ) \subseteq Y$ is open and contained in $L_j$, and we conclude that $\varphi(K_j^\circ) \subseteq L_j^\circ$. Then $K_j \subset K_{j+1}^\circ$ implies $L_j \subset L_{j+1}^\circ$ and the $L_j$ form indeed an exhaustion of compacts of $Y$.
\item Now $\varphi$ is indeed proper: A compact $L \subset Y$ is covered by the increasing $L_j^\circ, j \in \NN,$ and contained in one of them. \proven
\end{enumerate}
\end{proof}

\begin{proposition}
\label{prop_conjdisk}
Let $B$ be a complex space and $\varphi : \EE \to B$ a conjugating map. Then $\varphi$ is either biholomorphic or antibiholomorphic.
\end{proposition}
\begin{proof} \hfill
\begin{enumerate}
\item The group $\aut{\EE}$ acts transitively on $\EE$ and so does $\aut{B}$ on $B$. Therefore the structure around each point in $B$ is locally the same and $B$ is in fact a complex manifold. The Fatou-Bieberbach domains of $\EE$ form a sub-basis of the topology: For each $x_0 \in \EE$ and a ball $B(x_0, \delta) \subseteq \EE$ there is an $\varepsilon > 0$ such that $f(z) := \varepsilon \cdot z + x$ maps $\EE$ into $B(x, \delta)$. By lemma \ref{lem_homeo} $\varphi$ is a homeomorphism and therefore $B$ is a simply connected non-compact Riemann surface, and as such isomorphic to either $\EE$ or $\CC$. But a conjugating map $\varphi : \EE \to \CC$ cannot exist, because $\EE$ has a lot of Fatou-Bieberbach domains and $\CC$ has no Fatou-Bieberbach domain at all.
\item It is now sufficient to consider $\varphi : \EE \to \EE$ and we may assume $\varphi(0) = 0$ and $\varphi(1/2) = \alpha \cdot 1/2$ with $\alpha \in (0, 2) \subset \RR$. Recall that the automorphisms of $\EE$ are of the following form:
$$f_{a, \vartheta}(w) = e^{i \vartheta} \cdot \frac{a - w}{1 - \overline{a} w}, \quad a \in \EE, \; \vartheta \in [0, 2\pi)$$
An automorphism fixing $0$ is necessarly a rotation, therefore $\varphi$ maps rotations to rotations. The automorphisms can be decomposed into $f_{a, \vartheta} = f_{0, \vartheta} \circ f_{a, 0}$. The only elements of order two which are not rotations are the automorphisms $f_{a, 0}$, a fact also preserved under the conjugation by $\varphi$. Therefore the two types of automorphisms can be considered separately:
\begin{enumerate}
\item Rotations: $$\varphi(e^{i \vartheta} \cdot w) = e^{i \Theta(\vartheta)} \cdot \varphi(w), \quad w = \varphi^{-1}(z)$$
It follows that $e^{i \vartheta} \mapsto e^{i \Theta(\vartheta)}$ is multiplicative and continous. Therefore we have (up to ``gauging'' modulo $2\pi i \ZZ$) only two possibilities for $\Theta$, namely $\Theta(\vartheta) = \vartheta$ or $\Theta(\vartheta) = - \vartheta$.
\item Automorphisms $f_{a, 0}$: They are determined by the pre-image of $0$. Since $\varphi(0) = 0$, we immediately get that $f_{a, 0}$ is mapped to $f_{\varphi(a), 0}$.
\end{enumerate}
Putting this together we have:
$$\varphi\left( e^{i \vartheta} \cdot \frac{a - w}{1 - \overline{a} w} \right) = e^{ \pm i \vartheta} \cdot \frac{\varphi(a) - \varphi(w)}{1 - \overline{\varphi(a)} \varphi(w)}$$
In case of $e^{i \vartheta} \cdot w \mapsto e^{- i \vartheta} \cdot \varphi(w)$, we replace $\varphi$ by $\varphi \circ k$ where $k : \EE \to \EE$ denotes the complex conjugation. This changes
$\varphi\left( e^{i \vartheta} \cdot \frac{a - w}{1 - \overline{a} w} \right)$ to
$\varphi\left( e^{- i \vartheta} \cdot \frac{\overline{a} - w}{1 - a w} \right)$. We can now assume $e^{i \vartheta} \cdot w \mapsto e^{i \vartheta} \cdot \varphi(w)$.
In particular, this means that $\varphi$ preserves circles with orientation.
\item A calculation shows:
$$f_{a, 0} \circ f_{b, 0} (w) = \frac{a\overline{b} - 1}{1 - \overline{a}b} \cdot \frac{\frac{b - a}{1 - a \overline{b}} - w}{1 - \frac{\overline{b} - \overline{a}}{1 - \overline{a} b} w}$$
According to the previous step, $\varphi$ operates as follows:
$$\varphi\left(f_{a, 0} \circ f_{b, 0} (w)\right) = 
\frac{a\overline{b} - 1}{1 - \overline{a}b} \cdot \frac{\varphi(\frac{b - a}{1 - a \overline{b}}) - \varphi(w)}{1 - \overline{\varphi(\frac{b - a}{1 - a \overline{b}})} \varphi(w)} = 
\frac{a\overline{b} - 1}{1 - \overline{a}b} \cdot \frac{\frac{\varphi(b) - \varphi(a)}{1 - \overline{\varphi(b)}\varphi(a)} - \varphi(w)}{1 - \frac{\overline{\varphi(b)} - \overline{\varphi(a)}}{1 - \varphi(b)\overline{\varphi(a)}} \varphi(w)}$$
On the other hand, conjugation by $\varphi$ is a group homomomorphism:
$$\varphi\left(f_{a, 0} \circ f_{b, 0} (w)\right) = f_{\varphi(a), 0} \circ f_{\varphi(b), 0} (\varphi(w)) =
\frac{\varphi(a)\overline{\varphi(b)} - 1}{1 - \overline{\varphi(a)}\varphi(b)} \cdot \frac{\frac{\varphi(b) - \varphi(a)}{1 - \overline{\varphi(b)}\varphi(a)} - \varphi(w)}{1 - \frac{\overline{\varphi(b)} - \overline{\varphi(a)}}{1 - \varphi(b)\overline{\varphi(a)}} \varphi(w)}$$
Therefore it follows that
$$\frac{a\overline{b} - 1}{1 - \overline{a}b} = \frac{\varphi(a)\overline{\varphi(b)} - 1}{1 - \overline{\varphi(a)}\varphi(b)}$$
In case we choose $b \in \RR \cap \EE$ and set $a = 1/2$ this implies that $\varphi(b) \in \RR \cap \EE$, as $\varphi(1/2) \in \RR$ by assumption.
\item Because $\abs{\varphi}$ now only depends on the radius and the angle remains unchanged by $\varphi$, and because $\lim_{r \to 1} \varphi(r) = 1$, we can continously extend $\varphi$ to the boundary $\partial\EE$. The circle around $1/2$, going trough $0$ and $1$, is mapped to a circle around $1/2 \cdot \alpha$, still going through $0$ and $1$. Necessarly, it is $\alpha = 1$. By filling in circles of radius $2^{-k}, k \in \NN,$ with centers in $\RR \cap \EE$ we get a dense set of fixed points of $\varphi$ in $\RR \cap \EE$. By continuity and preserved unit scalar multiplication, $\varphi = \id_{\EE}$ follows, and $\id_{\EE} : \EE \to \EE$ is biholomorphic. \proven
\end{enumerate}
\end{proof}

\begin{theorem}
\label{thm_conjug3}
Let $X$ and $Y$ be complex spaces and $\varphi: X \to Y$ an iso-conjugating map. Then $\varphi$ is either biholomorphic or antibiholomorphic if the following criteria are fulfilled:
\begin{enumerate}
\item There is a holomorphic (not necessarly proper) embedding $X \hookrightarrow \EE^n$.
\item $X$ admits a holomorphic embedding $i : \EE \hookrightarrow X$ such that all holomorphic endomorphisms of $i(\EE) \isom \EE$ can be approximated uniformly on compacts by restrictions of holomorphic endomorphisms of $X$, and such that $i(\EE) = \bigcap_{k = 1, \dots, q} G_k^{-1}(x_0)$ for a point $x_0 \in X$ and $G_1, \dots, G_q \in \enmo{X}$.
\end{enumerate}
\end{theorem}
\begin{proof} The proof is similar to the one of theorem \ref{thm_conjug2}:
\begin{enumerate}
\item The image $B := \varphi(i(\EE))$ can be written as
\[\varphi(i(\EE)) = \bigcap_{k= 1, \dots, q} \left(\varphi \circ G_k \circ \varphi^{-1} \right)^{-1}(\varphi(x_0)),\] therefore it is an analytic subset of $Y$ (by lemma \ref{lem_analyticsub}), and carries the structure of a complex space.
\item We want to show that $\varphi|i(\EE) \to B$ is a conjugating map: By assumption the endomorphisms of the analytic subset $i(\EE) \isom \EE$ can be approximated by restrictions of holomorphic endomorphisms $F : X \to X$ and are conjugated to endomorphisms of $Y$ which restrict to endormorphisms of $B$. Because $\varphi$ is a homeomorphism, we also have approximation uniformly on compacts for the conjugated endomorphisms, and therefore $\varphi|i(\EE)$ is a conjugating map. Then by proposition \ref{prop_conjdisk} $\varphi|i(\EE)$ is either biholomorphic or antibiholomorphic to its image.
\item The algebras $\holalgb{X}$ and $\holalgb{Y}$ of bounded holomorphic functions are isomorphic as $\CC$-algebras by sending $\holalgb{X} \ni f \mapsto (2\norm{f}) \cdot (\varphi \circ i)^{-1} \circ (\varphi \circ i \circ f/(2\norm{f}) \circ \varphi^{-1}) = f \circ \varphi^{-1} =: \hat{\Phi}(f) \in \holalgb{Y}$, if necessary with complex conjugation. The properties of a $\CC$-algebra homomorphism are fulfilled, since addition and multiplication of functions is point-wise. $\hat{\Phi}$ is an isomorphism because the situation in $X$ and $Y$ is symmetric and because $\varphi$ is an iso-conjugating map. In case that an additional complex conjugation is necessary, we replace $X$ by its image under complex conjugation in $\EE^n \subset \CC^n$.
\[
\begin{diagram}
\node[2]{X} \arrow{e,t}{\varphi} \arrow{sw,t}{f/(2\norm{f})} \arrow{s,r}{} \node{Y} \arrow{s,l}{} \arrow{se,t}{f/(2\norm{f}) \circ \varphi^{-1}} \\
\node{\EE} \arrow{e,t}{i} \node{X} \arrow{e,t}{\varphi} \node{Y} \node{\EE} \arrow{w,t}{\varphi \circ i}
\end{diagram}
\]
\item By assumption, there is a holomorphic embedding $(f_1, \dots, f_n) : X \hookrightarrow \EE^n$ with $f_1, \dots, f_n \in \holalgb{X}$. It follows that $f_k \circ \varphi^{-1} \in \holalgb{Y}, \; k=1,\dots,n$ and that $(f_1 \circ \varphi^{-1}, \dots, f_n \circ \varphi^{-1}) : Y \hookrightarrow  \EE^n$ is a holomorphic embedding too (but not neccesarly proper). As the map $(f_1, \dots, f_n)$ is biholomorphic to its image, this implies that $\varphi^{-1}$ is holomorphic. Now we apply the same argument to $\hat{\Phi}^{-1} : \holalgb{Y} \to \holalgb{X}$ and obtain that $\varphi$ is holomorphic too. \proven
\end{enumerate}
\end{proof}

\begin{examples} \hfill
\begin{enumerate}
\item $\EE^n, \; n \in \NN$
\item $\EE \times \Omega$, where $\Omega \subset \CC^m$ is a bounded open set: $i : \EE \hookrightarrow \EE \times \Omega, \; z \mapsto (z, x_0), \; x_0 \in \Omega$ is a holomorphic embedding such that $i(\EE) = f^{-1}(0, x_0)$ where $f(z, w) = (z, x_0)$ is a holomorphic endomorphism of $\EE \times \Omega$. Any $\alpha : i(\EE) \to i(\EE)$ is the restriction of $A : \EE \times \Omega \to \EE \times \Omega, \; A(z, w) = (\alpha(z), w)$.
\item The unit ball $B_n(0,1) \subset \EE^n$: $i : \EE \hookrightarrow B_n(0, 1), \; z \mapsto (z, 0, \dots, 0)$ is a holomorphic embedding such that $i(\EE) = f^{-1}(0,\dots,0)$ where $f(z_1, \dots, z_n) = (z_1, 0, \dots, 0)$ is a holomorphic endomorphism of $B_n(0, 1)$. Any $\alpha : i(\EE) \to i(\EE)$ is the restriction of $A : B_n(0, 1) \to B_n(0, 1), \; A(z_1, z_2, \dots, z_n) = (\alpha(z), 0,\dots,0)$.
\item $\EE^n \! \setminus \! \{ 0 \}, \; n \geq 2$: $i : \EE \hookrightarrow \EE^n \! \setminus \! \{ 0 \}, \; z \mapsto (z, 1, \dots, 1)$ is a holomorphic embedding such that $i(\EE) = f^{-1}(1,\dots,1)$ where $f(z_1, \dots, z_n) = (1, z_2, \dots, z_n)$ is a holomorphic endomorphism of $\EE^n \! \setminus \! \{ 0 \}$. Any $\alpha : i(\EE) \to i(\EE)$ is the restriction of $A : \EE^n \! \setminus \! \{ 0 \} \to \EE^n \! \setminus \! \{ 0 \}, \; A(z_1, z_2, \dots, z_n) = (\alpha(z), 1,\dots,1)$.
\end{enumerate}
\end{examples}

\section{Open questions}
\label{sec_questions}
\subsection{Smallest sub-semigroups charaterizing a bounded domain}
\label{sec_questions_smallest}
We note that in general much smaller sub-semigroups of endomorphisms are sufficient to consider. For example, let $\Omega_1, \Omega_2 \subset \CC$ be two bounded domains and $\varphi : \Omega_1 \to \Omega_2$ an iso-conjugating map. This map is a homeomorphism and induces a homeomorphism $\Phi : \EE \to \EE$ between the universal coverings $p : \EE \to \Omega_1$ and $q : \EE \to \Omega_2$. The endomorphisms $f : \Omega_1 \to \Omega_1$ also induce endomorphisms $F : \EE \to \EE$ such that the following diagram commutes:
\[
\begin{diagram}
\node[3]{\EE} \arrow[6]{e,t}{\Phi} \arrow[3]{s,l}{p} \arrow[2]{sw,t}{F} \node[6]{\EE} \arrow[3]{s,r}{q} \arrow[2]{sw,t}{\Phi \circ F \circ \Phi^{-1}} \\ \\
\node{\EE} \arrow[6]{e,t}{\Phi} \arrow[3]{s,l}{p} \node[6]{\EE} \arrow[3]{s,r}{q} \\
\node[3]{\Omega_1} \arrow[6]{e,t}{\varphi} \arrow[2]{sw,t}{f} \node[6]{\Omega_2} \arrow[2]{sw,r}{\varphi \circ f \circ \varphi^{-1}} \\ \\
\node{\Omega_1} \arrow[6]{e,t}{\varphi} \node[6]{\Omega_2}
\end{diagram}
\]
Therefore the map $\Phi : \EE \to \EE$ is a conjugating map for some sub-semigroup of $\enmo{\EE}$. From Eremenko's theorem we know that $\varphi : \Omega_1 \to \Omega_2$ is biholomorphic, and we can conclude that $\Phi$ is biholomorphic too.
This raises the question how small sub-semigroups of the endomorphism semigroup need to be in order to guarantee (anti-)holomorphicity. For the unit disk for example, it is enough to consider certain contractions for a countable dense set of points in order to get continuity for the conjugating map and after that, the automorphism group is large enough to show the (anti-)holomorphicity.
\subsection{Stein manifolds with density property}
\label{sec_questions_density}
Varolin \cite{Varolin1} introduced the notion of the density property for complex manifolds:
\begin{definition}
A complex manifold $X$ has the \textit{density property}, if the Lie algebra generated by the completely integrable holomorphic vector fields on $X$ is dense in the Lie algebra of all holomorphic vector fields on $X$, where dense is meant with respect to the compact-open topology. 
\end{definition}
The idea is somehow to ensure that such a complex manifold with density property has ``a lot'' of automorphisms. Additionally, there will also be a lot of endomorphisms, as there are many holomorphic embeddings $\CC^n \hookrightarrow X$ (see proposition \ref{prop_localcontraction} below), and a lot of functions on $X$, provided that it is Stein too. Of course, such manifolds are good candidates for being characterized by their endomorphism semigroup.
\par 
The examples \ref{ex_conjug1}.\ref{ex_semisimple} (semi-simple homogenous spaces $X=G/K$ where $G$ is a semi-simple Lie group with trivial center and $K$ is a reductive subgroup), \ref{ex_conjug1}.\ref{ex_alggroup} (linear algebraic groups) and \ref{ex_conjug1}.\ref{ex_hypersurface} (certain hypersurfaces) represent the known classes of examples for Stein manifolds with density property (see Toth and Varolin \cite{Varolin3}, Kaliman and Kutzschebauch \cite{KalimanKutzschebauch1}, \cite{KalimanKutzschebauch2}). This leads to the following
\begin{question}
Do all manifolds with density property admit a proper holomorphic embedding of the complex line?
\end{question}
If true, all such manifolds would be determined by their endomorphism semigroup.
Another question in this context is:
\begin{question}
Do the Fatou-Bieberbach domains of a Stein manifold with density property form a sub-basis of the topology?
\end{question}
This is only known for $\CC^n, n \geq 2,$ as a consequence of the previously mentioned result of Buzzard and Hubbard \cite{Buzzard1}. This would also imply that these manifolds are determined by their endomorphism semigroup, because of the following:
Varolin \cite{Varolin2} showed that for each point in a manifold with density property, there is an automorphism with this point as an attractive fixed point. In addition, he generalized the following proposition of Rosay and Rudin \cite{RosayRudin} (originally for $X = \CC^n$):
\begin{proposition}
\label{prop_localcontraction}
Let $X$ be a complex manifold, $\kappa \in \mathrm{Aut}\left(X\right)$, $x_0 \in X$ with $\kappa(x_0) = x_0$ and the eigenvalues $\lambda_i$ of $d_{x_0}\kappa$ satisfy $\abs{\lambda_i} < 1$. Then
\[
U := \left\{x \in X \,:\, \lim_{r \to \infty} \kappa^r(x) = x_0 \right\}
\]
is a domain, biholomorphic to $\CC^n$.
\end{proposition}
Together, this results in: $\forall \, x_0 \in X \; \exists \, i : \CC^n \hookrightarrow U \subseteq X$ holomorphic embedding with open image and $x_0 = i(0)$. As a domain of attraction of an automorphism, $U \subseteq X$ is Runge. Therefore all holomorphic endomorphisms of $U$ can be approximated by endomorphisms of $X$, considered as holomorphic functions (via the embedding $i$). An iso-conjugating map $\varphi : X \to Y$ can now be assumed to be a homeomorphism, because $Y$ is a manifold (since $\aut{X} \isom \aut{Y}$ act transitively -- see Varolin \cite{Varolin2}) and lemma \ref{lem_homeo} can be applied. Therefore, also $\varphi(U) \subset Y$ is Runge and it follows that $\enmo{U} \isom \enmo{\varphi(U)}$. By the already established result (example \ref{ex_conjug1}.\ref{ex_cn}) for $\CC^n \isom U$, we can conclude that $\varphi|U$ is either biholomorphic or antibiholomorphic. Thus, $\varphi$ itself is either biholomorphic or antibiholomorphic on each connected component.
\section{Acknowledgement}
I would like to thank Frank Kutzschebauch for helpful and interesting discussions and critical remarks.
\bibliographystyle{amsplain}

\end{document}